\documentclass[a4paper,10pt]{article}
\usepackage[utf8]{inputenc}

\usepackage{amsmath,amssymb,enumerate}
\usepackage{graphicx}
\usepackage{caption}
\usepackage{subcaption}

\newcommand{\R}{\mathbb{R}}
\newcommand{\pd}{\partial}

\renewcommand{\Re}{\text{Re}}
\newcommand{\Wh}{\text{Wh}}
\newcommand{\ordo}{\mathcal{O}}

\newcommand{\sep}{; }

\newtheorem{thm}{Theorem}

\newtheorem{dfn}{Definition}

\begin{document}

\title{
Numerical reconstruction of pulsatile blood flow
from 4D computer tomography angiography data
}

\author{
Attila~Lovas$^1$\footnote{lovas@math.bme.hu}, 
R\'obert~Nagy$^2$, 
Elek~Csobo$^1$,
\\
Brigitta~Szil\'agyi$^3$,
P\'eter~S\'otonyi$^4$
}

\maketitle

\begin{center}
${}^1$Budapest University of Technology and Economics, 
Department of Analysis, Building H,
Egry~J\'ozsef str. 1, Budapest, Hungary 1111

\medskip
${}^2$Budapest University of Technology and Economics,
Department of Structural Mechanics, Building K, 
M\H{u}egyetem rkp. 3, Budapest, Hungary 1111

\medskip
${}^3$Budapest University of Technology and Economics, 
Department of Geometry, Building H, 
Egry~J\'ozsef str. 1, Budapest, Hungary 1111

\medskip
${}^4$Semmelweis University, Heart and Vascular Center,
Department of Vascular Surgery, V\'arosmajor str. 68,
Budapest, Hungary 1122 
\end{center}

\begin{abstract}

We present a novel numerical algorithm developed to 
reconstuct pulsatile blood flow from 
ECG-gated CT angiography data.
A block-based optimization method was constructed
to solve the inverse problem corresponding to the
Riccati-type ordinary differential equation that
can be deduced from conservation principles and 
Hooke's law. Local flow rate for $5$ patients
was computed in $10$ cm long
aorta segments that are located $1$ cm below the heart.
The wave form of the local flow rate curves seems to be 
realistic.
Our approach is suitable for estimating 
characteristics of pulsatile blood flow in aorta based on 
ECG gated CT scan thereby 
contributing to more accurate description of
several cardiovascular lesions.
\footnote{Submitted to Computers and Mathematics With Applications, 2015. 
(lovas@math.bme.hu)
}


\end{abstract}

\section{Introduction}

Pulsatile flow in blood vessels has been studied for more than
300 years. Leonard Euler initiated the theory of pressure wave
propagation in the vascular system in 1775. The first modern mathematical
model of pulsatile flow in blood vessels was developed by 
Diederik Johannes Korteweg and Horace Lamb.
Hitherto, numerous attempts have been made to study blood flow in
the vascular system. Yunlong Huo and Ghassan S. Kassab have published
a Womersley-type mathematical model to analyze pulsatile blood flow in
the entire coronary arterial tree \cite{YunlongHuo2006}. 
Mostly stenosed geomeries were investigated owing to their
connection to vascular diseases such as
arteriosclerosis which may lead to heart attack and stroke
\cite{MoloyKumar2012}. Mir Golam Rabby et al.
have investigated numerically the effects of non-Newtonian
modeling on unsteady periodic flows in a two-dimensional
pipe with idealized stenoses \cite{Golam2014}.
In the above mentioned studies, motions of the vessel wall
were not taken into consideration.

ECG-gated computer tomography along with an
image segmentation algorithm based on Active Contour Model
\cite{Terzo1988}
enables us to measure the local cross sectional area of the vessel
as a function of time and the position of the considered 
cross section along the vessel's axis. 
The main objective of the present work is to develop
a numerical procedure to compute volumetric flow rate 
of blood through the local cross sections of the vessel 
from the local cross sectional area.
We assume that there exist a short segment of the vessel where
it is thin and elastic. Furthermore, the change in the 
local cross sectional area must be caused by the change 
in blood pressure only, and  gravitation has no effect 
on the fluid because the patient is
lying during the medical examination.
By applying the one-dimensional Euler momentum equation,
equation of continuity and Hooke's law, we
deduced a Riccati-type ODE for volumetric flow rate at the
beginnig of the considered vessel segment closer to the heart.
The only periodic solution with positive integral was accepted as
physically admissible. Existence of such a solution can be tested
numerically. However, the Riccati-type ODE
contains an extra scalar parameter and can be solved just for
fixed parameter values. To find the volumetric flow rate a
block-based optimization method was implemented.

\section{Material and methods}

\subsection{ECG-gated CT angiography}


Imaging of a pulsatile organ is a highly demanding application for any cross 
sectional imaging modality. Computed tomography (CT) imaging of the heart 
became widely available with the introduction of multi-detector CT (MDCT) 
scanners with four-slice detector arrays and 500 ms minimum rotation time 
\cite{CRBecker2000}. Still images of the moving heart was generated using 
retrospective ECG-gating: slow table motion during spiral scanning and 
simultaneous acquisition of the slices and the digital ECG trace provided 
oversampling of scan projections \cite{CRBecker2000}.
After the exposure, slices 
recorded in the same phase of the ECG trace are matched to generate a 
3D dataset of the volume of interest, representing either systole or diastole. 
A drawback of this method is the higher radiation dose compared to the 
normal non-oversampled spiral acquisition \cite{JPEarls2008}. 
An important advantage is the possibility to reconstruct multiphasic 
datasets of the same volume, resulting in motion images of the same slices. 
This allows us the functional analysis of the moving organs such as the heart 
or the large vessels. Recent advances in CT technology 
($256-320$ detector rows, $270$ ms minimum rotation time) allow 
for rapid ECG-gated CTA of the whole aorta during a single breath-hold. 

Imaging of the aorta was performed in $5$ patients ($5$ men, mean age $68.2\pm 6.1$ years)
with a $256$-slice MDCT (Philips Brilliance iCT, Koninklijke Philips N.V., Best,
The Netherlands) using a retrospectively ECG-gated protocol tailored for the
imaging of the aorta. Investigations were performed on images readily available from
patients with suspected aortic disease. Low-dose (tube voltage: $100$ kV) native
scan was followed by a retrospective ECG-gated CT angiography of the whole 
aorta ($100$ kV) with a reduced field of view to maximize spatial resolution. 
Nonionic contrast agent was injected into an antecubital vein at a flow rate
of $4-5$ ml/s using a power injector. Images were reconstructed using a sharp
convolution kernel and iterative reconstruction algorithm (iDose4, Koninklijke
Philips N.V., Best, The Netherlands) with a slice thickness of $1$ mm and an 
increment of $1$ mm. Multiphase images were reconstructed corresponding every 
10\% of the R-R cycle resulting in ten series of images for each patient. 
Patients gave written informed consent before CT examination was performed.
Experimental protocol and informed consent was approved by the 
Regional Ethical Committee of Semmelweis University ($133/2011$).

\subsection{Image segmentation}

The first problem which we have to solve in order to compute
cross sectional area of the blood vessel is to locate it in
the dicom image. This can be carried out by some sort of image segmentation
algorithm. In our case, Actice Contour Model (ACM)
was applied to perform segmentation. The first version
of ACM was published by Michael Kass, Andrew Witkin
and Demetri Terzopoulos \cite{Terzo1988} at the end of the 80s.

Active contours are energy minimizing splines in 2D
which they are commonly referred to as snakes.
It is possible to generalize the algorithm to higher
dimensions using energy minimizing surfaces but
the corresponding numerical method is usually
not stable because a surface might intersect itself
or change its topology. For that reason, an algorithm based 
on a sequential two-dimensional balloon snake model
was applied to each phase which means that
the final position of the moving snake on a
slice was implemented as initial condition on the next one.
As the result of the segmentation process, we get
the coordinates of the vessel's internal wall
for each slice perpendicular
to $z$ axis in the region of interest and for all phases in time.
For more information about the Active Contours we refer the reader to
\cite{Sonka2006} and \cite{Hlavac2007}.

To mitigate the spatial measurement error, we fitted a bicubic smoothing 
spline to the point cloud resulting from the active contour method in each 
time step. This way the increase of temporal resolution also became possible 
exploiting the affine covariance of B--splines and ensuring the periodic 
movement of the control points by trigonometric approximation using 
the first 3 harmonics of the heart rate \cite{Nagy2015}. 
The cross-sectional area is defined as the area of the section 
perpendicular to the center line of the fitted surface.

\subsection{Governing Equations}

To make our exposition self-contained we present some laws 
of continuum mechanics, these are the conservation laws for 
mass and momentum. In addition, empirical constitutive laws are needed 
to relate certain unknown variables such as relations between stress and strain. 
 
Let us examine the pulsating flow of blood in an artery assuming that 
the wall is thin and elastic. Owing to the pressure gradient the artery
wall deforms and the elastic restoring force of the wall makes it possible 
for waves to propagate and so it maintains a pulsating motion of the artery.
The artery radius $r(t,x)$ varies in time along the artery. 
Let the local cross sectional area be $S=\pi r^2$ and the averaged
flow velocity $u(t,x)$ are periodic in the first
argument with period $T$ where $T$ denotes the period of the cardiac
cycle. We define the volume flow rate as $q=u S$. 
Assuming that blood is homogeneous and incompressible we obtain from the 
law for conversation of mass that

\begin{equation}\label{eq:1}
 \frac{\pd S}{\pd t} + \frac{\pd q}{\pd x} = 0.
\end{equation}
Additionally the law for conservation of momentum can be written as follows:
\begin{equation}\label{eq:2}
 \frac{\pd u}{\pd t} + u \frac{\pd u}{\pd x}
 =
 -\frac{1}{\rho}\frac{\pd p}{\pd x}.
\end{equation}

We use Hooke's law to relate stress and strain rates. 
Let $h$ be the vessel wall thickness, assumed to be much smaller 
than the vessel radius, and Young's modulus will be denoted by $E$. 
The change in tube radius must be caused by the blood pressure. 
The elastic strain due to the lengthening of the circumference is $dr/r$. 
The change is elastic force must be balanced by the changing in 
pressure force $2a dp$, which implies

\[
\frac{dp}{dr}={Eh}{r^2} \quad \text{and} \quad \frac{dp}{dS}=\frac{\sqrt{\pi}Eh}{2 S^{3/2}}.
\]
From this we obtain the following:
\begin{equation}\label{eq:3}
 \frac{\pd p}{\pd x}
 =
 \frac{hE\sqrt{\pi}}{2 S^{3/2}}\frac{\pd S}{\pd x}.
\end{equation}

\section{Calculation}

\subsection{Periodic Riccati equation}

Let us consider a nonbranching cylindral blood
vessel segment of length $L$ where equations
\eqref{eq:1}, \eqref{eq:2} and \eqref{eq:3} are satisfied.
If we substitute \eqref{eq:3} into \eqref{eq:2} we have
\begin{equation}\label{eq:4}
 \frac{\pd u}{\pd t} + u \frac{\pd u}{\pd x}
 =
 -\frac{\alpha}{S\sqrt{S}}\frac{\pd S}{\pd x},
\end{equation}
where
$
\alpha = 0.5hE\sqrt{\pi}/\rho
$
is assumed to be constant during the cardiac cycle. 

If we multiply the left-hand side of \eqref{eq:4}
by $S$, using $q = uS$ and integrating by parts
we obtain
\begin{equation}
\begin{split}
\int\limits_0^L S\frac{\pd u}{\pd t} + q \frac{\pd u}{\pd x}
\,\mathrm{d}x 
=&
\left[\frac{q^2(t,x)}{S(t,x)}\right]_{x=0}^{x=L} 
+
\int\limits_0^L S\frac{\pd u}{\pd t}
+ u\frac{\pd S}{\pd t}
\,\mathrm{d}x
\\
=&
\left[
\frac{q^2(t,x)}{S(t,x)}
\right]_{x=0}^{x=L}
+
\int\limits_0^L \frac{\pd q}{\pd t}
\,\mathrm{d}x.
\end{split}
\end{equation}
By applying  \eqref{eq:1}, the volumetric flow rate 
can be written as
\begin{equation}
 q(t,x) = Q(t) + \Phi(t,x),
\end{equation}
where
$Q(t)=q(t,0)$
and
$$
\Phi (t,x) = - \int\limits_0^x \frac{\pd S}{\pd t} (t,y)
\,
\mathrm{d}y.
$$
From this we obtain the following Riccati-type
ordinary differential equation
\begin{equation}\label{eq:ric}
 \frac{\mathrm{d} Q}{\mathrm{d} t} = A Q^2 + B Q + C,
\end{equation}
where coefficients are periodic with period $T$ and can be
written as follows:
\begin{align}
 A(t) &= 
 -\frac{1}{L}
 \left[
 \frac{1}{S(t,x)}
 \right]_{x=0}^{x=L}
 \\
 B(t) &=
 -\frac{2}{L}
 \frac{\Phi(t,L)}{S(t,L)}
 \\
 C(t) &=
 -\frac{1}{L}
 \left(
 \frac{\Phi(t,L)^2}{S(t,L)}
 -
 2\alpha\left[
 \sqrt{S(t,x)}
 \right]_{x=0}^{x=L}
 +
 \int\limits_0^L
 \frac{\pd \Phi}{\pd t}
 \,\mathrm{d}x
 \right).
\end{align}

\subsection{Periodic solutions}

\paragraph{Physically admissible solutions}

Assume that $x=0$ corresponds to the beginnig of the considered vessel segment
closer to the heart.
\begin{dfn}
 A solution of \eqref{eq:ric} $Q$ is said to be
 \emph{physically admissible} if it satisfies the following
 conditions:
 \begin{enumerate}[i.]
  \item $Q:\R \to \R$ is a periodic function with period $T$,
  \item $\int\limits_0^T Q(t)\,\mathrm{d}t>0$.
 \end{enumerate}

\end{dfn}

It is well known that the general solution of a scalar Riccati equation can
be obtained by quadrature whenever at least one particular solution $Q_0$ is
known. The substitution
\begin{equation}
 Q = Q_0 - \frac{1}{W}
\end{equation}
in the Riccati equation yields
\begin{equation}
 \frac{\mathrm{d} W}{\mathrm{d} t} = -(2Q_0 A + B)W + A,
\end{equation}
which is first order inhomogeneous linear ODE and its
general solution can be written as
\begin{equation}
 W(t) = K W_h(t) + W_{ih} (t),
\end{equation}
where $W_h$ is the solution of the homogeneous equation for which
$W_h(0)=1$ holds and $W_{ih}$ is 
the solution of the inhomogeneous equation
for which
$W_{ih}(0)=0$ holds. Periodicity of $Q$ requires that
$
Q(0) = Q(T)
$
that holds if and only if $K$ solves the
quadratic equation
\begin{equation}\label{eq:14}
 Q_0(T) K^2 + 
 \left(
 Q_0 (T) \frac{W_{ih}(T)}{W_h(T)}
 -
 \frac{W_h(T)-1}{W_h(T)}
 \right) K
 - \frac{W_{ih}(T)}{W_h(T)} = 0,
\end{equation}
where $Q_0(0)$ is supposed to be zero.
The number of periodic solutions depends on the sign of the
discriminant that can be tested numerically.
\begin{equation}
 \Delta = \left(
 Q_0 (T) \frac{W_{ih}(T)}{W_h(T)}
 -
 \frac{W_h(T)-1}{W_h(T)}
 \right)^2 + 4 Q_0(T) \frac{W_{ih}(T)}{W_h(T)}
\end{equation}
Although the non-negativity of $\Delta$ guarantees the existence of
periodic solution, it does not imply that at least one of the
solutions is physically admissible.


Leon Kotin demonstrated the existence and uniqueness of a
positive and of a negative periodic solution of a periodic Riccati
equation in which the coefficients satisfy certain general condition.
Moreover, any solution which is everywhere continuous lies between these 
two solutions, and every solution is asymptotic to one of these as 
the independent variable increases or decreases \cite{Kotin19681227}.
More precisely the following is true.
\begin{thm}[Kotin 1968] 
 If $A',B,C$ in \eqref{eq:ric} are continuous everywhere and $-AC>0$
 holds, then equation \eqref{eq:ric} has a unique positive
 solution $Q_+$ and a unique negative solution $Q_{-}$ which are 
 periodic with period $T$ and any solution behaves asymptotocally like
 $Q_+$ or $Q_-$ as $t\to\infty$.
\end{thm}

\paragraph{Assymptotic stability}
It is clear that if conditions of
Kotin's theorem are satisfied, then $Q_+$ is the unique 
physically admissible solution. Hovewer, $Q_+$ may not be
assimptotically stable. Consider the case when $C<0$ and $Q$ is an
arbitrary perturbation of $Q_+$. From the unicity of
the positive solution, we have that there exists $t_0\in\R$
where $Q$ vanishes and $Q$ can have only one root,
since at $Q=0$ the right-hand side of \eqref{eq:ric} is equal to
$C$ which is negative.
As a consequence, we have that the
physically admissible solution is not necessarily 
assimptotically stable.

\subsection{Numerical Procedure}

\paragraph{The block-based opmimization method}

For any fixed values of $\alpha$, the functions
$W_h$ and $W_{ih}$ are computed numerically.
In cases where $\Delta\ge 0$, periodic solutions of \eqref{eq:ric}
can be
expressed as follows
$$
Q_{1,2}(t) = Q_0(t) - \frac{1}{K_{1,2}W_h(t)+W_{ih}(t)},
$$
where $K_1$ and $K_2$ are solutions of \eqref{eq:14}. 
In our case, $\alpha$ is not known thus the problem is underdetermined up to
a constant. To find the local flow rate further conditions
are needed.

To overcome this problem, we consider a second non-branching vessel segment of length $L$ that is
adjacent to the previous one such that the end
of the first segment is joined to the beginning of the second one. 
Quantities corresponding to the second
vessel segment are denoted by tilded symbols
(e.g. $\widetilde{q}$ and $\widetilde{Q}$) and $\alpha$ is
incorporated as a subscript to indicate the dependence of the
considered quantity on the model parameter.
We define the so-called "internal consistency" functional
\begin{equation}\label{eq:func}
 I(\alpha) = \int\limits_0^T (\widetilde{q}_\alpha(t,0)-q_\alpha (t,L))^2 \,\mathrm{d}t
\end{equation}
that measures the mean squared difference between the local flow rate at
the beginning of the second vessel segment computed from 
motions of the second vessel segment and
the local flow rate at
the end of the first vessel segment 
calculated from data corresponding to the first vessel segment.
The domain of $I$ is 
$$
\text{dom}(I)=
\left\lbrace \alpha\in\R^+ 
\left| 
\Delta_\alpha\ge 0
\,
\&
\,
\widetilde{\Delta}_\alpha
\ge 0
\right.  
\right\rbrace.
$$
Let us define the average flow rate
corresponding to $\alpha$ by
$$
\bar{q}_{\alpha}=
\frac{1}{2T} 
\int\limits_0^T
\widetilde{q}_\alpha (t,0)
+
q_\alpha(t,L)
\,\mathrm{d}t.
$$
Our goal is to minimalize $I(\alpha)$ for $\alpha$-s
that satisfy
$
\bar{q}_{min}
\le
\bar{q}_{\alpha}
\le
\bar{q}_{max}
$,
where $q_{min}$ and $q_{max}$
are the minimal and maximal average flow rate
which may occur
under physiological circumstances.

\paragraph{Sensitivity analysis}

At this point we examine the effect of small 
perturbations of $\alpha$ around the optimal value
($\alpha_{opt}$). We can write
\begin{equation}
 C = C_0 + \alpha C_1,
\end{equation}
where $C_0$ and $C_1$ do not depend on $\alpha$.
First-order approximation of $Q_{\alpha}$ can be
written as
\begin{equation}
 Q_\alpha (t) = Q_{\alpha_{opt}}(t) + P(t)(\alpha-\alpha_{opt})
 +
 \ordo
 \left(|\alpha-\alpha_{opt}|^2\right),
\end{equation}
where $P$ is called sensitivity coefficient 
and it is the unique periodic solution 
of the following 
first-order linear differential
equation:
\begin{equation}
 \frac{\mathrm{d}P}{\mathrm{d}t}
 = (2AQ_{\alpha_{opt}}+B)P+C_1.
\end{equation}

\section{Results}

\subsection{Local flow rates}

Gender, age, body mass index (BMI), heart rate (HR)
and dose length product (DLP) are presented in Table \ref{tab:1}.
Calculations were performed on $10$ cm long non-branching 
segments of the descending aorta.
Considered vessel segments are located $1$ cm below the heart.
In each case, the last $2\times 1$ cm long parts of the vessel segment
provides input data for the block-based optimization algorithm.
According to the physiological fact that the cardiac output
lies between 
$4$ $\text{dm}^3/\text{min}$ and 
$6$ $\text{dm}^3/\text{min}$, 
maximal average flow rate ($\bar{q}_{max}$) 
was set to $100$ $\text{cm}^3/\text{s}$ and
the minimal average flow rate ($\bar{q}_{min}$) 
to $66.7$ $\text{cm}^3/\text{s}$.
MATLAB's ode45 function was used to solve differential
equations that may occur in simulations.
\begin{table}[!h]
 \centering
 \begin{tabular}{c|c|c|c|c|c}
        No. & Gender  & Age  & BMI    & HR   & DLP   \\
  \hline 1. & M       & $69$ & $29.1$ & $60$ & $3504$\\
  \hline 2. & M       & $67$ & $23.9$ & $52$ & $2690$\\
  \hline 3. & M       & $62$ & $28.9$ & $86$ & $2756$\\
  \hline 4. & M       & $65$ & $19.3$ & $59$ & $2161$\\
  \hline 5. & M       & $78$ & $24.7$ & $67$ & $2300$\\
 \end{tabular}
 \caption{
 Patient data.
 Gender -- male (M)/female (F), Age (years), 
 BMI ($\text{kg}/\text{m}^{2}$), 
 HR (bpm), DLP (mGycm).
 }
 \label{tab:1}
\end{table}

MATLAB's fzero function was applied to compute $\alpha$ values corresponding
to the minimal and maximal flow rate.
These are $\alpha_{\text{min}}$ and $\alpha_{\text{max}}$ respectively.
We found that iterations converged to a limit in all cases.
Global minimum of $I(\alpha)$ was computed on the interval
$(\alpha_{\text{min}},\alpha_{\text{max}})$ by MATLAB's
fminbnd function. 
We found that interations converged to a minimum in all cases
and the interval $(\alpha_{\text{min}},\alpha_{\text{max}})$ is a
subset of the domain of 
$I$ that is the feasible set of the optimization problem.
Values of $\alpha_{\text{min}}$, $\alpha_{\text{max}}$
and $\alpha_{\text{opt}}$ in $10^3\text{cm}^3/\text{s}^2$ 
are presented in Table \ref{tab:2}.
The mean squared error between blocks  
$\text{MSE}=(I(\alpha_{\text{opt}})/T)^{1/2}$
is introduced to measure the goodness of the model (Table \ref{tab:2}).
Conditions of Kotin's theorem were tested and the
results are presented in Table \ref{tab:2}.
\begin{table}[!h]
\centering
 \begin{tabular}{c|c|c|c|c|c}
  No. & $\alpha_{\text{min}}$ & $\alpha_{\text{max}}$ & $\alpha_{\text{opt}}$ & MSE  & Kotin\\
 \hline  1. & $2.72$ & $7.56$ & $2.75$  & $2.25$ & + \\
 \hline  2. & $2.47$ & $6.87$ & $2.47$  & $3.42$ & + \\
 \hline  3. & $2.67$ & $7.42$ & $2.68$  & $3.82$ & + \\
 \hline  4. & $2.81$ & $7.99$ & $2.88$  & $3.06$ & + \\
 \hline  5. & $5.81$ & $15.8$ & $6.82$  & $2.12$ & - \\
 \end{tabular}
 \caption{Parameter values corresponding to 
 the minimal, maximal and the optimal local flow rate 
 ($10^3 \text{cm}^3/\text{s}^2$), MSE ($\text{cm}^3/\text{s}$),
 and conditions of Kotin's theorem -- satisfied (+)/ not satisfied (-).
 }
 \label{tab:2}
\end{table}

\begin{figure}
    \centering
    \begin{subfigure}[!h]{0.9\textwidth}
        \includegraphics[width=\textwidth]{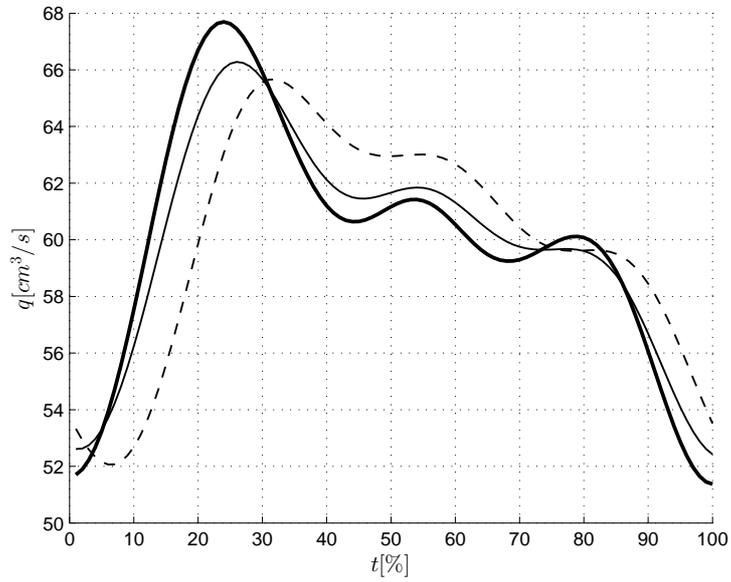}
        \caption{Patient No. 3}
        \label{fig:qYoung}
    \end{subfigure}
    ~ 
    \begin{subfigure}[!h]{0.9\textwidth}
        \includegraphics[width=\textwidth]{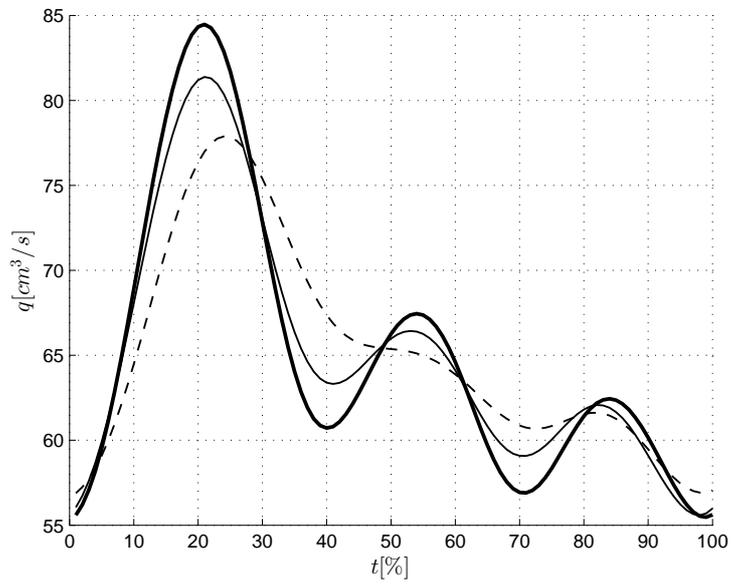}
        \caption{Patient No. 5}
        \label{fig:qOld}
    \end{subfigure}   
    \caption{
    Local flow rate ($\text{cm}^3/\text{s}$) in time 
    (\% of the cardiac cycle). 
    Thick  line -- 10\%,  
    solid  line -- 50\%
    dashed line -- 90\%
    of the aorta segment.
    }
    \label{fig:q}
\end{figure}

Detailed results are presented for the youngest 
(patient No. $3$) and for the oldest patient (patient No. $5$).
Local flow rates were calculated at the $10\%$, $50\%$ and
$90\%$ of the aorta segment.
We plot the local flow rates corresponding to the youngest and the oldest
patient at different equally distributed times during the cardiac cycle in
Fig. \ref{fig:qYoung} and Fig. \ref{fig:qOld}).

\begin{figure}
    \centering
    \begin{subfigure}[!h]{0.9\textwidth}
        \includegraphics[width=\textwidth]{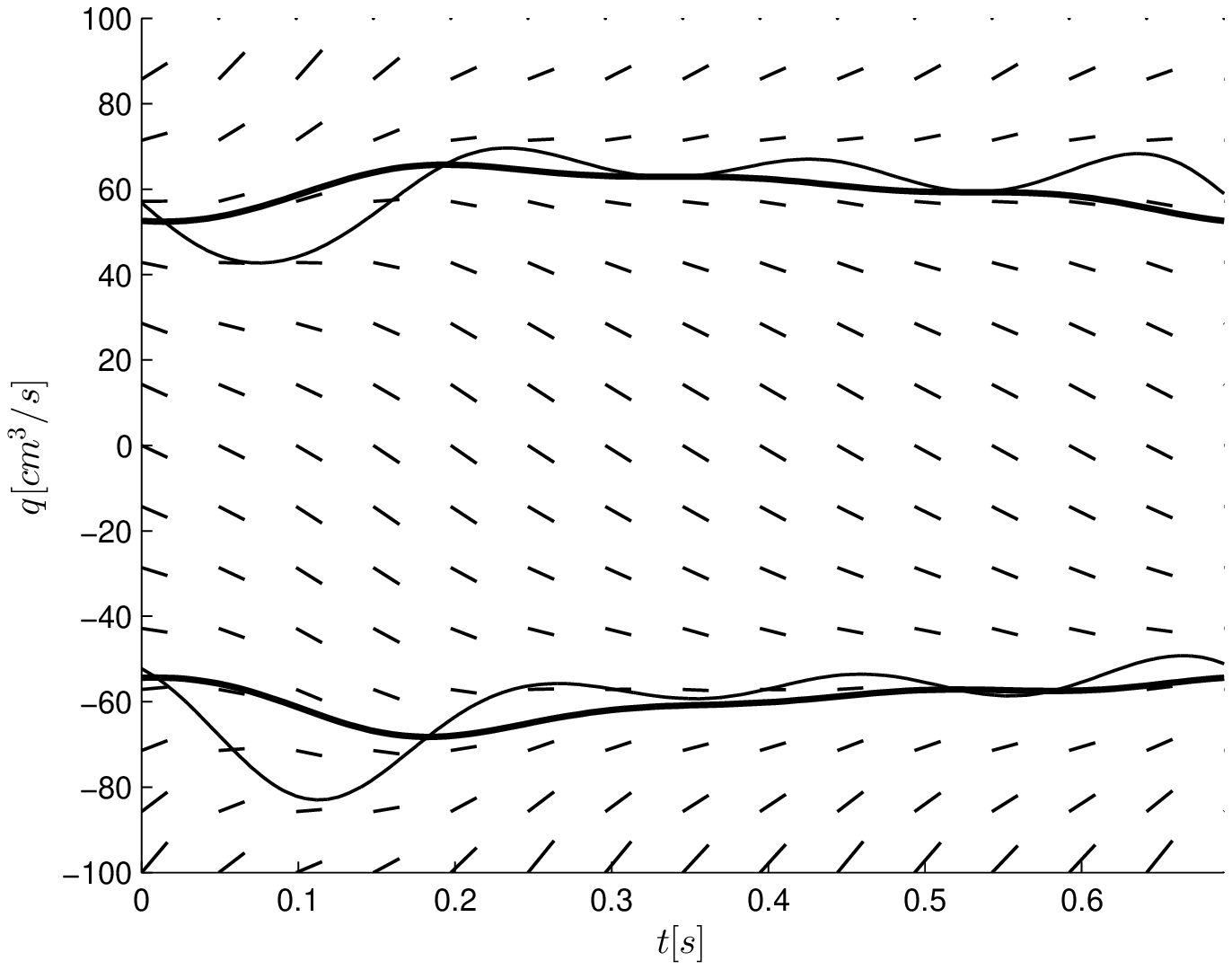}
        \caption{Patient No. 3}
        \label{fig:phYoung}
    \end{subfigure}
    ~ 
    \begin{subfigure}[!h]{0.9\textwidth}
        \includegraphics[width=\textwidth]{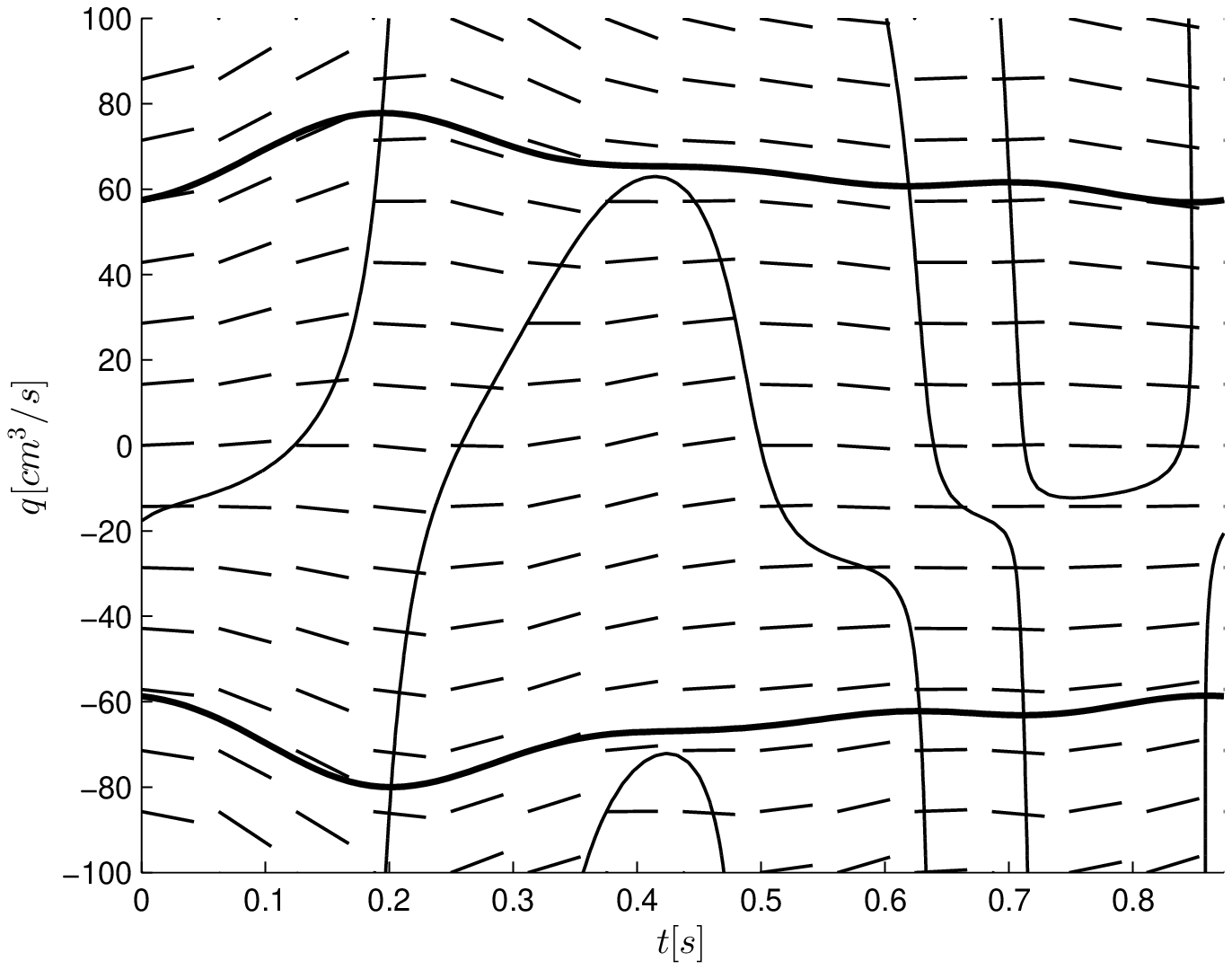}
        \caption{Patient No. 5}
        \label{fig:phOld}
    \end{subfigure}   
    \caption{
    Phase pictures. Local flow rate ($\text{cm}^3/\text{s}$) 
    in time (s). 
    Thick line -- periodic solution,
    solid line -- nullcline.
    }
    \label{fig:ph}
\end{figure}

\smallskip
We have seen that conditions of Kotin's theorem are satisfied 
in almost all cases and the only exception is just patient No. 5
who was the oldest person. The typical behaviour of the solutions
can be seen in Fig. \ref{fig:phYoung} and the phase picture 
corresponding to the only non-typical case is illustrated in
Fig. \ref{fig:phOld}.

\begin{figure}
    \centering
    \begin{subfigure}[h]{0.9\textwidth}
        \includegraphics[width=\textwidth]{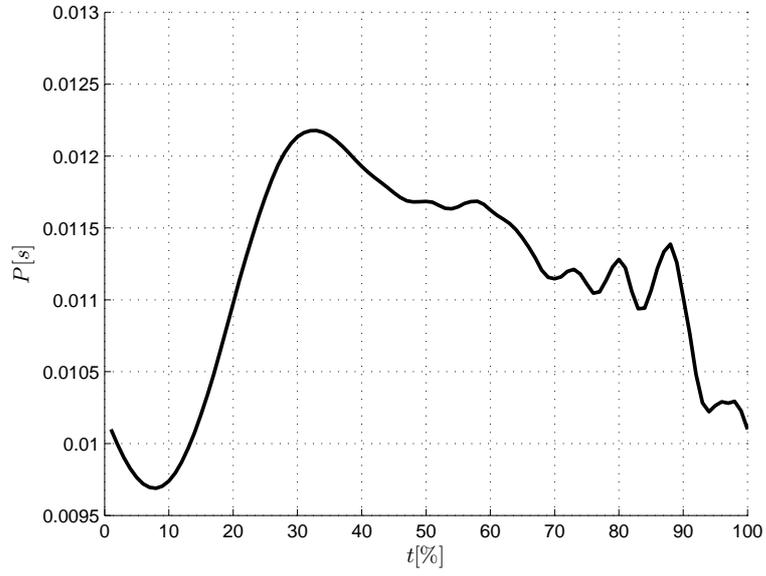}
        \caption{Patient No. 3}
        \label{fig:sYoung}
    \end{subfigure}
    ~ 
    \begin{subfigure}[h]{0.9\textwidth}
        \includegraphics[width=\textwidth]{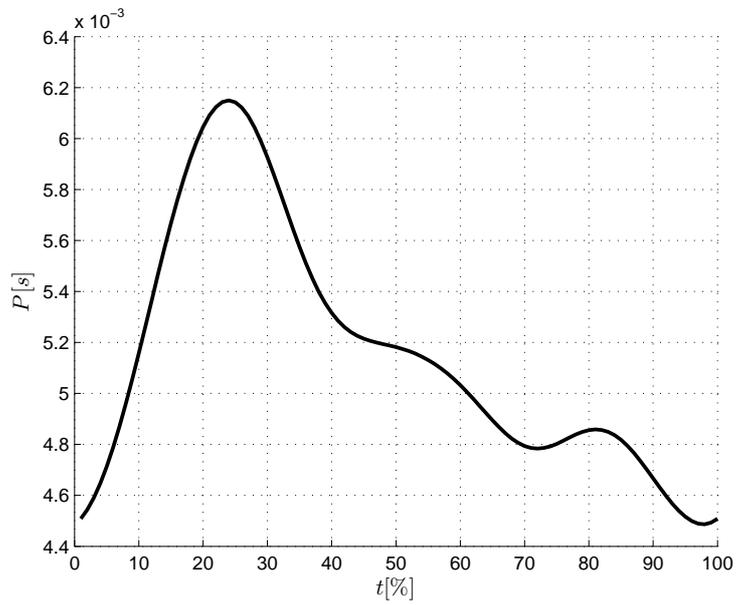}
        \caption{Patient No. 5}
        \label{fig:sOld}
    \end{subfigure}   
    \caption{
    Sensitivity coefficient 
    $P$ (s) in time (\% of the cardiac cycle).
    }
    \label{fig:s}
\end{figure}

The graph of the sensitivity coefficients can be seen in
Fig. \ref{fig:sYoung} and Fig. \ref{fig:sOld}.

\subsection{Reynolds and Womersley numbers}

\begin{figure}
  \centering
    \begin{subfigure}[!h]{0.9\textwidth}
        \includegraphics[width=\textwidth]{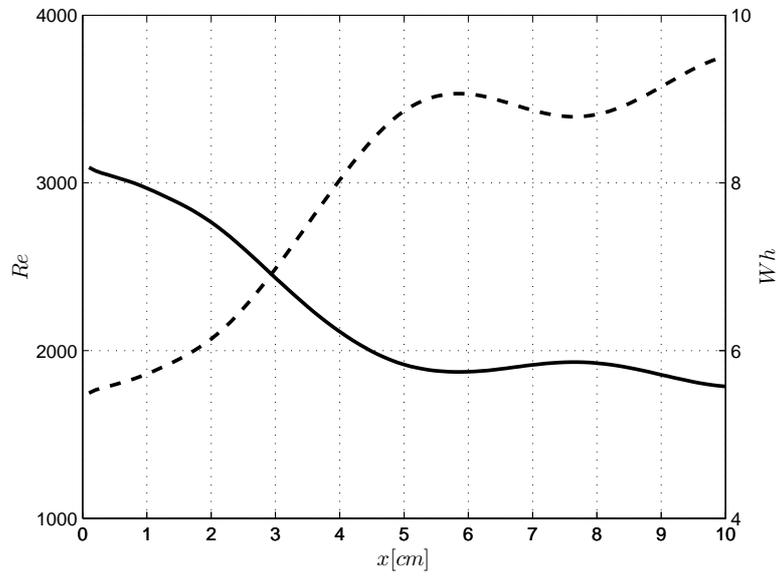}
        \caption{Patient No. 3}
        \label{fig:wrYoung}
    \end{subfigure}
    \begin{subfigure}[!h]{0.9\textwidth}
        \includegraphics[width=\textwidth]{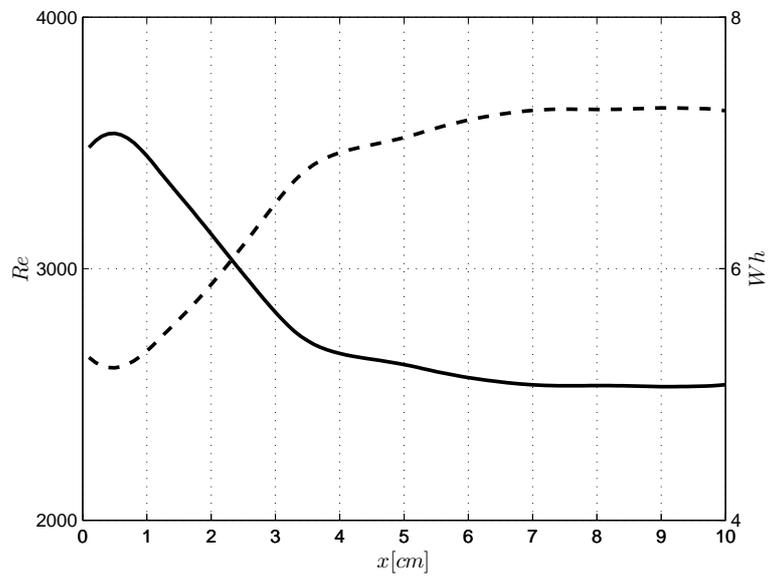}
        \caption{Patient No. 5}
        \label{fig:wrOld}
    \end{subfigure}  
    \caption{
    Reynolds and Womersley numbers
    vs. distance from the beginning of 
    the aorta segment (cm).
    }
    \label{fig:wr}
\end{figure}

For the sake of completeness, we give the definitions of Reynolds ($\Re$) and
Womersley ($\Wh$)
numbers. The first is a dimensionless number in fluid mechanics
that is defined as the ratio of inertia forces to viscous forces,
expressed in tubular flows as
\begin{equation}
 \Re = \frac{2 v r \rho}{\eta},
\end{equation}
where $v$ is the flow velocity, $r$ is the vessel's radius, $\rho$
is the blood density and $\eta$ is the dynamical viscosity of blood.
The second is a dimensionless number in biofluid mechanics
that relates the transient inertial forces to viscous effects. 
The Womersley number is defined by
\begin{equation}
 \Wh = 2r\sqrt{\frac{\omega \rho}{\eta}},
\end{equation}
where $\omega$ is the angular frequency of the oscillations. 
In accordance with the literature 
\cite{Jung2014}, \cite{Hinghofer1987}, 
in our calculations
blood density 
was set to $1.06$ $\text{g}/\text{cm}^3$
and the
kinematical viscosity of blood
to $3.5\times 10^{-3}$ $\text{Pa}\cdot\text{s}$.
Estimatons for Reynolds and Womersley numbers along the artery segment 
are illustrated in common coordinate system 
in Fig. \ref{fig:wrYoung} and Fig. \ref{fig:wrOld}.

The Womersley number is a dynamic similarity measure of oscillatory flows 
relating inertia and viscous forces. In rigid pipes for laminar incompressible 
flows small values (approx. Wo$<$1) allow the development of the parabolic 
velocity profile of the steady state solution and the flow is almost in phase 
with the pressure gradient, while large values (approx. Wo$>$10) indicate a 
flat velocity profile with a good approximation and the flow follows 
the pressure gradient by about 90 degrees in phase. 
In the presented examples -- as seen in Fig.~\ref{fig:wrYoung} 
and \ref{fig:wrOld} -- the values lie inbetween, 
yielding a complex time-dependent velocity profile.

The Reynolds number is another dynamic similarity measure relating inertia 
and viscous forces. In rigid pipes small values (approx. Re$<$2100) indicate 
laminar flow, while large values (approx. Re$>$4000) correspond to turbulent 
flow. In the transition zone, where also our example is situated, the behavior
strongly depends on the existing disturbances in the flow.

\section{Discussion}

Our work present a novel numerical algorithm for the non-invasive
blood flow quantification based on ECG-gated CT angiography images.

Quantification of flow in blood vessels provides important information 
aiding the management of different clinical conditions \cite{Powell2000}. 
Since its introduction in the late 1980s, phase-contrast 
magnetic resonance (MR) imaging has become the method of choice 
for flow quantification over ultrasound, as it is able to measure 
time-resolved 3D flow with high accuracy and reproducibility without 
the need of a good acoustic window \cite{Stankovic2014}.
Our approach enables us to 
calculate aortic flow on a routine ECG-gated CT angiography dataset. 
This is of huge clinical potential, as the knowledge of hemodynamic 
parameters could vastly improve the diagnostic performance of CT imaging 
in several cardiovascular pathologies, such as aortic coarctation or dissection.  
Our method can also be beneficial in the management of patients
who are not candidates for an MR examination with long acquisition 
time either because of claustrophobia or due to unstable patient status. 

State-of-the-art investigations of the disorders of human arterial sections 
apply 3D numerical 
fluid--structure interaction simulations involving the calculation of 
the blood flow field inside the lumen, the necessary boundary conditions 
of which are the time-dependent pressure profile at the outlet and volumetric
flow rate at the inlet cross-section. The protocol to determine 
these functions non-invasively is the measurement of both on the arm 
followed by transforming them to the section under scrutiny by a 
1D systemal model of the circulatory system. 
Our method offers a different and easy to use procedure to formulate 
the inlet boundary condition without Dopler velocimetry, 
thus facilitating 3D simulations.

\section{Conclusions}

We presented a novel numerical algorithm for the non-invasive 
quantification of blood flow based on ECG-gated CT angiography images. 
This method can vastly improve the diagnostic performance of CT imaging
in several cardiovascular pathologies. Further studies are needed 
to validate our approach in a clinical setting. 

A block-based optimization algorithm was developed to
compute local flow rate in time.
The wave form of the local flow rate curves
seems to be
similar 
to those 
that can be measured by invasive intraarterial method.
The mean squared error (MSE) values,
which measures the goodness of the model,
are quite small compared to the average flow rate.

According to the analysis of Reynolds and Womersley numbers,
the velocity profile has a complex time-dependent behaviour in
the detailed examples (Patient No. 3 and No. 5).

We have found that conditions of Kotin's theorem are satisfied in most
of the cases and the physically admissible solution is not 
assymptotically stable in general. Further mathematical analysis of 
equation \eqref{eq:ric} is needed to find weaker sufficient
conditions for the physical admissibility.

\section*{References}

\bibliography{pulsatility}

\end{document}